\documentclass[12pt,a4paper]{article}
\usepackage{amscd,amsmath,amssymb,amsfonts,times}
\usepackage{mathrsfs}
\usepackage[all]{xy}
\numberwithin{equation}{subsection}
    \newtheorem{thm}{Theorem}[subsection]
    \newtheorem{lem}[thm]{Lemma}
    
    \newtheorem{prop}[thm]{Proposition}
    \newtheorem{cor}[thm]{Corollary}
    
    \newtheorem{claim}[thm]{Claim}

    \newtheorem{rem}[thm]{Remark}

\newcommand{\qed}
{\mbox{}\nolinebreak$\square$\medbreak\par}
\newenvironment{pf}{\par\smallskip\noindent\emph{Proof.}}{\hfill\qed\par\smallskip}
\newenvironment{pf*}[1]{\par\smallskip\noindent\emph{#1.}}{\hfill\qed\par\smallskip}
\begin{document}
\title{Quintic surface over $p$-adic local fields with 
infinite $p$-primary torsion in the Chow group of $0$-cycles}
\author{Masanori Asakura}
\date\empty
\maketitle

\def\can{\text{\rm can}}
\def\codim{\text{\rm codim}}
\def\cano{\text{\rm canonical}}
\def\cd{\text{\rm cd}}
\def\ch{\text{\rm ch}}
\def\Coker{\mathrm{Coker}}
\def\Cor{\text{\rm Cor}}
\def\crys{\mathrm{crys}}
\def\dlog{d{\text{\rm log}}}
\def\dR{\text{\rm d\hspace{-0.2pt}R}}
\def\Eis{{\text{\rm Eis}}_{\text{\rm form}}}
\def\et{\text{\rm \'et}}  % etale
\def\exp{\text{\rm exp}}
\def\Frac{\text{\rm Frac}}
\def\phami{\phantom{-}}
\def\id{\text{\rm id}}              % identity
\def\Image{\text{\rm Im}}        % image
\def\Hom{\text{\rm Hom}}  
\def\Mor{\text{\rm Mor}}  
\def\ker{\text{\rm Ker}}          % kernel
\def\mf{\text{\rm mapping fiber of}}
\def\Pic{\text{\rm Pic}}
\def\CH{\text{\rm CH}}
\def\NS{\text{\rm NS}}
\def\End{\text{\rm End}}
\def\pr{\text{\rm pr}}
\def\Proj{\text{\rm Proj}}
\def\ord{\text{\rm ord}}
\def\qis{\text{\rm qis}}
\def\reg{\text{\rm reg}}          %
\def\res{\text{\rm res}}          %
\def\Res{\text{\rm Res}}
\def\Spec{\text{\rm Spec}}     % spectrum
\def\syn{\text{\rm syn}}
\def\tors{\text{\rm tors}}
\def\tr{\text{\rm tr}}
\def\cont{\text{\rm cont}}
\def\zar{\text{\rm zar}}
\def\AJ{\text{\rm AJ}}
\def\Div{\mathrm{div}}
\def\an{\mathrm{an}}
\def\bA{{\mathbb A}}
\def\bC{{\mathbb C}}
\def\C{{\mathbb C}}
\def\G{{\mathbb G}}
\def\bE{{\mathbb E}}
\def\bF{{\mathbb F}}
\def\F{{\mathbb F}}
\def\bH{{\mathbb H}}
\def\bJ{{\mathbb J}}
\def\bN{{\mathbb N}}
\def\bP{{\mathbb P}}
\def\P{{\mathbb P}}
\def\bQ{{\mathbb Q}}
\def\Q{{\mathbb Q}}
\def\bR{{\mathbb R}}
\def\R{{\mathbb R}}
\def\bZ{{\mathbb Z}}
\def\Z{{\mathbb Z}}

\def\cC{{\mathscr C}}
\def\cD{{\mathscr D}}
\def\cE{{\mathscr E}}
\def\cK{{\mathscr K}}
\def\cO{{\mathscr O}}
\def\O{{\mathscr O}}
\def\cR{{\mathscr R}}
\def\cS{{\mathscr S}}
\def\cU{{\mathscr U}}
\def\cV{{\mathscr V}}
\def\cX{\hspace{-0.8pt}{\mathscr X}\hspace{-1.2pt}}

\def\fo{\mathfrak o}

\def\XR{\cX}
\def\XK{X}
\def\XZp{\cX}
\def\XQp{X}
\def\XbK{\ol X}

\def\UR{\cU}
\def\UK{U}
\def\UZp{\cU}
\def\UQp{U}
\def\UbK{\ol U}
\def\UF{\UR_s}

\def\SR{\cV}
\def\SK{S}
\def\SbK{\ol S}
\def\SF{\SR_s}

\def\CR{\cC}
\def\CK{C}
\def\CZp{\cC}
\def\CQp{C}
\def\CbK{\ol C}

\def\DR{\cD}
\def\DK{D}
\def\DZp{\cD}
\def\DQp{D}
\def\DbK{\ol D}
\def\tDZp{\widetilde{\cD}}
\def\XbK{\ol X}

%                                 Greece
%
\def\ep{\epsilon}
\def\vG{\varGamma}
\def\vg{\varGamma}
\def\M{\mathscr{M}}
\def\uC{\mathscr{C}}
\def\uG{\mathscr{G}}
\def\uJ{\mathscr{J}}

%
%
%
%
%
%                                 simple
%
%
%
\def\lra{\longrightarrow}
\def\lla{\longleftarrow}
\def\Lra{\Longrightarrow}
\def\hra{\hookrightarrow}
\def\lmt{\longmapsto}
\def\ot{\otimes}
\def\op{\oplus}
%                              decolation
\def\wt#1{\widetilde{#1}}
\def\wh#1{\widehat{#1}}
\def\spt{\sptilde}
\def\ol#1{\overline{#1}}
\def\ul#1{\underline{#1}}
\def\us#1#2{\underset{#1}{#2}}
\def\os#1#2{\overset{#1}{#2}}
\def\lim#1{\us{#1}{\varinjlim}}
\def\plim#1{\us{#1}{\varprojlim}}
%
%                         extra-ordinary symbol
%
%
%
\def\Gm{{\mathbb G}_{\hspace{-1pt}\mathrm{m}}}
\def\Ga{{\mathbb G}_{\hspace{-1pt}\mathrm{a}}}
\def\zp{\bZ_p}
\def\qp{\bQ_p}
\def\qzp{\qp/\zp}
\def\ql{\bQ_{\ell}}
\def\zl{\bZ_{\ell}}
\def\qzl{\ql/\zl}
\def\qz'{{\bQ}/{\bZ}'}
\def\isom{\hspace{9pt}{}^\sim\hspace{-16.5pt}\lra}
\def\lisom{\hspace{10pt}{}^\sim\hspace{-17.5pt}\lla}
%{\hspace{9pt}\sptilde\hspace{-16.5pt}\lra}

\def\hA{\wh{A}}
\def\uA{\mathscr A_0}
\def\uuA{\mathscr A}
\def\uB{\mathscr B_0}
\def\uuB{{\mathscr B}}
\def\uK{{{\mathscr R}^\flat_K}}
\def\uoK{{{\mathscr R}^\flat_{\ol K}}}
\def\uuK{{\uuR_K}}
\def\uR{{{\mathscr R}^\flat}}
\def\uuR{{\mathscr R}}
\def\tD{\widetilde{D}}
\def\tT{\widetilde{T}}
\def\tZ{\widetilde{Z}}

\def\psr{R\hspace{.5pt}[\hspace{-1.57pt}[q_0]\hspace{-1.58pt}]}
\def\psro{R\hspace{.5pt}[\hspace{-1.57pt}[q_1]\hspace{-1.58pt}]}
\def\psf{R\hspace{.5pt}(\hspace{-1.57pt}(q_0)\hspace{-1.58pt})}
\def\psfo{R\hspace{.5pt}(\hspace{-1.57pt}(q_1)\hspace{-1.58pt})}
\def\psru{\hspace{.5pt}[\hspace{-1.57pt}[u]\hspace{-1.58pt}]}
\def\psrt{\hspace{.5pt}[\hspace{-1.57pt}[t]\hspace{-1.58pt}]}
\def\psrq{\hspace{.5pt}[\hspace{-1.57pt}[q_0]\hspace{-1.58pt}]}
\def\psfu{\hspace{.5pt}(\hspace{-1.57pt}(u)\hspace{-1.58pt})}
\def\psfq{\hspace{.5pt}(\hspace{-1.57pt}(q_0)\hspace{-1.58pt})}
\def\psrqq{\hspace{.5pt}[\hspace{-1.57pt}[q]\hspace{-1.58pt}]}
\def\psrqqq{\hspace{.5pt}[\hspace{-1.57pt}[q_0,q_0^{-1}]\hspace{-1.58pt}]}
\def\psfqq{\hspace{.5pt}(\hspace{-1.57pt}(q)\hspace{-1.58pt})}
\def\psrqu{\hspace{.5pt}[\hspace{-1.57pt}[q_0,u]\hspace{-1.58pt}]}
\def\psfti{\hspace{.5pt}(\hspace{-1.57pt}(t_i)\hspace{-1.58pt})}
\def\psftt{\hspace{.5pt}(\hspace{-1.57pt}(t-1)\hspace{-1.58pt})}
\def\psfqi{\hspace{.5pt}(\hspace{-1.57pt}(q_i)\hspace{-1.58pt})}

\def\tet{{\tau^\et_\infty}}
\def\tdR{{\tau^\dR_\infty}}
\def\tsyn{{\tau^\syn_\infty}}
\def\htsyn{\wh{\tau}_{\infty}^{\,\syn}}
\def\htet{\wh{\tau}_{\infty}^{\,\et}}
\def\ttsyn{\wt{\tau}_{\infty}^{\syn}}
\def\ttet{\wt{\tau}_{\infty}^{\et}}
%
%
%
%                             Hodge-Witt sheaf
%
\def\mwitt#1#2#3{W_{\hspace{-2pt}#2}{\hspace{1pt}}\omega_{#1}^{#3}}
\def\witt#1#2#3{W_{\hspace{-2pt}#2}{\hspace{1pt}}\Omega_{#1}^{#3}}
\def\mlogwitt#1#2#3{W_{\hspace{-2pt}#2}{\hspace{1pt}}\omega_{{#1},{\log}}^{#3}}
\def\tlogwitt#1#2#3{W_{\hspace{-2pt}#2}{\hspace{1pt}}\tom_{{#1},{\log}}^{#3}}
\def\logwitt#1#2#3{W_{\hspace{-2pt}#2}{\hspace{1pt}}\Omega_{{#1},{\log}}^{#3}}
\def\loglogwitt#1#2{W_{\hspace{-2pt}n}{\hspace{1pt}}\Omega^{#2}_{Y^{(#1)}\hspace{-1.5pt},\hspace{1pt}E^{[#1]}\hspace{-1pt},\hspace{1pt}\log}}

\section{Introduction}
In the paper \cite{AS} S. Saito and the author
constructed a surface $X$ over a $p$-adic local field such that the $l$-primary torsion
part $\CH_0(X)[l^\infty]$ of the Chow group of $0$-cycles
is infinite for $l\ne p$, which gave a counter-example
to a folklore conjecture.
The purpose of this paper is to show that there is such an example
even for $l=p$:
\begin{thm}[Theorem \ref{main00}]\label{main0}
There is a quintic surface $X\subset \P^3_{\Q_p}$ over $\Q_p$ such that the $p$-primary torsion part
$\CH_0(X\times_{\Q_p}K)[p^\infty]$ is not finite for arbitrary finite extension $K$ of $\Q_p$.
\end{thm}
Our proof is comparable with that of \cite{AS}, however
a new difficulty appears in case $l=p$.
Let us recall the outline of the proof of \cite{AS} briefly.
It follows from the universal coefficient theorem on Bloch's
higher Chow group that we have the exact sequence
\begin{equation}\label{bloch}
0\lra \CH^2(X,1)\ot\Q_l/\Z_l\os{i}{\lra}
\CH^2(X,1;\Q_l/\Z_l)\lra 
\CH^2(X)[l^\infty]\lra0
\end{equation}
for any $l$ (possibly $l=p$).
The proof of \cite{AS} breaks up into two steps.
We first showed that 
if $X$ is generic then
$\CH^2(X,1)\ot\Q_l/\Z_l$ contains only decomposable elements
supported on hyperplane section
(cf. Lem. \ref{as1} below).
Next we showed that the {\it boundary map}
$$
\partial:\CH^2(X,1;\Q_l/\Z_l)\lra \Pic(Y)\ot\Q_l/\Z_l
$$
is surjective (modulo finite groups) for $X$ which has a good reduction $Y$.
Thus if $Y$ contains primitive divisors, then 
$\CH^2(X,1;\Q_l/\Z_l)$ contains indecomposable elements and hence
the map $i$ cannot be surjective.

The technique used in the former step works also in case $l=p$.
On the other hand, in the latter step, we used the result of Sato-Saito \cite{SS},
in which they proved a weak Mordell-Weil type theorem for Chow group mod $l$
different from $p$.
Unfortunately its mod $p$ counterpart has not been obtained. 
Thus we cannot use the same technique as in \cite{AS} 
to show the surjectivity of
the boundary map in case $l=p$.

Actually we do not need the surjectivity of $\partial$ to prove Theorem \ref{main0}.
It is enough to show that the corank of the image of $\partial$ is greater than one.
To do this, we construct an indecomposable element in
$\CH^2(X,1;\Z/p^n\Z)$ (never coming from $\CH^2(X,1)\ot\Z/p^n\Z$ !). 
The strategy is as follows. 
We consider a quintic surface $X\subset \P^3_{\Q_p}$ which contains
an irreducible quintic curve $C$ with four nodes. 
Let $\tilde{C}\to C$ be the normalization and
$\{P_i,Q_i\}$ $(1\leq i \leq 4)$ the inverse images of the four nodes on $C$.
Since the curve $\tilde{C}$ is genus 2, there is a rational function $f_n$ 
on $\tilde{C}$ such that
$\Div(f_n)\equiv\sum_{i=1}^4 r_i(P_i-Q_i)$ mod $p^n$
(this is a simple application of the theorem of Mattuck \cite{mattuck}
which asserts that the Jacobian $J(\tilde{C})(\Q_p)$ is isomorphic to
$\Z_p^2$ modulo finite groups). Thus the pair $(C,f_n)$ determines an element
in $\CH^2(X,1;\Z/p^n\Z)$. We then prove that its boundary is nontrivial (hence
indecomposable) under some assumptions.

\medskip

This paper is organized as follows. In \S \ref{nodesect}, we give an axiomatic
approach to the construction of surface with infinite $p$-primary torsion
in the Chow group of $0$-cycles.
In \S \ref{quinticsect},
we will construct such a quintic surface over $\Q_p$. 
A technical difficulty appears in the calculations
of the boundary map. There we will use Igusa's $j$-invariants of
hyperelliptic curves of genus two, which we list in Appendix for the convenience
of the reader.

\section{Preliminaries}
For an abelian group $M$ we denote by $M[n]$ (resp. $M/n$) 
the kernel (resp. cokernel) of the multiplication by $n$.
We denote the $p$-primary torsion by $M[p^\infty]=\cup_{n\geq 1} M[p^n]$.
For schemes $X$ and $T$ over a base scheme $S$, we write $X(T)=\Mor_S(T,X)$
the set of $S$-morphisms,
and say $x\in X(T)$ a $T$-valued point of $X$.
If $T=\Spec R$, then we also write $X(R)=X(\Spec R)$ and say $x\in X(R)$
a $R$-rational point.
 
For a regular scheme $X$,
we denote by $Z_i(X)=Z^{\dim X-i}(X)$ the free abelian group of irreducible subvarieties
of Krull dimension $i$.
\subsection{$K$-cohomology and Gersten complex}
Let $X$ be a smooth variety over a field $F$.
Let us denote by $X^i$
the set of irreducible subvarieties of $X$ of codimension $i$.
We write the function field of $Z$ by $\eta_Z$.

Let $\cK_i$ be the sheaf associated to a presheaf $U\mapsto K_i(U)$
where $K_i(U)$ is Quillen's $K$-theory.
The Zariski cohomology group $H^\bullet(X,\cK_i)$ is called the $K$-cohomology.
We only concern with $H^1(X,\cK_2)$. It has an explicit description 
by using the Gersten complex
\begin{equation}\label{gersten0}
\begin{CD}
K^M_2(\eta_X)
@>{d_2}>>
\bigoplus_{\codim ~D=1} \eta_D^\times
@>{d_1}>>
Z^2(X).
\end{CD}
\end{equation}
Recall the maps $d_1$ and $d_2$.
We denote by $(f,D)$ the image of an element $f\in \eta_D^\times$ via the 
canonical inclusion $\eta_D^\times\to\bigoplus_{\codim ~D=1} \eta_D^\times$.
Then the map $d_2$ (called the {\it tame symbol}) is defined as follows
$$
d_2\{f,g\}=\sum_{\codim ~D=1}
\left((-1)^{\ord_D(f)\ord_D(g)}\frac{f^{\ord_D(g)}}{g^{\ord_D(f)}}|_D,D\right).
$$
The map $d_1$ is defined in the following way.
Let $\tilde{D}\to D$ be the normalization.
To $f\in \eta_D^\times$, we associate the Weil divisor $\Div_{\tilde{D}}(f)$
on $\tilde{D}$. Letting $j:\tilde{D}\to D\hra X$ be the composition,
$d_1(f,D)$ is defined to be $j_*(\Div_{\tilde{D}}(f))$.
It is simple to check $d_1d_2=0$.
Tensoring \eqref{gersten0} with $\Z/n$, one has a complex
\begin{equation}\label{gersten2}
\begin{CD}
K^M_2(\eta_X)/n
@>{d_2\ot\Z/n}>>
\bigoplus_{\codim ~D=1} \eta_D^\times/n
@>{d_1\ot\Z/n}>>
Z^2(X)/n
\end{CD}
\end{equation}
Then the Gersten conjecture (Quillen's theorem) tells that 
\eqref{gersten2} gives rise to a flasque resolution of the shaef $\cK_2/n$
and hence one has the canonical
isomorphism
\begin{equation}\label{gersten3}
H^1(X,\cK_2/n)\cong \ker(d_1\ot\Z/n)/\Image(d_2\ot\Z/n)
\end{equation}
for each $n\geq 0$.
Hereafter we always identify the $K$-cohomolgy $H^1(X,\cK_2/n)$
with the group in the right hand side of \eqref{gersten3}.
\subsection{$\CH^2(X,1)$ and $K$-cohomology}
We denote by $\CH^i(X,j;G)$
Bloch's higher Chow group with coefficients in an abelian group $G$.
We simply write $\CH^i(X,j)=\CH^i(X,j;\Z)$ and $\CH^i(X)=\CH^i(X,0)$. 
By \cite{lands} 2.5, we have the canonical isomorphism
\begin{equation}\label{lands}
\CH^2(X,1;\Z/n)\cong H^1(X,\cK_2/n)
\end{equation}
for each $n\geq0$. We will also identify 
$\CH^2(X,1;\Z/n)$ with $H^1(X,\cK_2/n)$ by the above isomorphism.

By the universal coefficients theorem on higher Chow group
there is the exact sequence
\begin{equation}\label{bloch0}
0\lra \CH^2(X,1)/n
\lra
\CH^2(X,1;\Z/n)\lra 
\CH^2(X)[n]\lra0
\end{equation}
for $n\ne0$.
Putting $n=l^k$ and taking the inductive limit on $k$, one obtains \eqref{bloch}.
Suppose that $n$ is prime to the characteristic of $F$.
Then there is the regulator map
$$
\mathrm{reg}_X:\CH^2(X,1;\Z/n)\lra H^3_\et(X,\Z/n(2))
$$
to the etale cohomology group.
Let
\begin{align*}
NH^3_\et(X,\Z/n(2)):&=\ker(H^3_\et(X,\Z/n(2))\lra
H^3_\et(\eta_X,\Z/n(2)))\\
&=\Image(\bigoplus_{\codim~Z=1}H^3_Z(X,\Z/n(2))\lra
H^3_\et(X,\Z/n(2)))
\end{align*}
where the second equality follows from the localization exact sequence of
etale cohomology.
The main theorem of Bloch-Ogus theory tells that
the image of the regulator map coincides with the above:
\begin{equation}\label{bloch1}
\CH^2(X,1;\Z/n)\os{\sim}{\lra} NH^3_\et(X,\Z/n(2)).
\end{equation}
It follows from \eqref{bloch0} and \eqref{bloch1} that one has
Bloch's exact sequence
\begin{equation}\label{bloch2}
0\lra \CH^2(X,1)/n
\lra
NH^3_\et(X,\Z/n(2))
\lra 
\CH^2(X)[n]\lra0.
\end{equation}
If $F$ is a $p$-adic local field, the cohomology group 
$H^\bullet_\et(X,\Z/n(j))$ (resp. $H^\bullet_\et(X,\Q_l/\Z_l(j))$)
is known to be finite (resp. of cofinite type).
Hence the $n$-torsion part $\CH^2(X)[n]$ is finite and
$\CH^2(X)[l^\infty]$ is of cofinite type for any $n$ and $l$:
$$
\CH^2(X)[l^\infty]\cong (\Q_l/\Z_l)^{\op r}+\mbox{(finite group)}.
$$

\subsection{Boundary map}\label{boundarysect}
Let $R$ be a discrete valuation ring with a prime element $\pi$.
Put $K:=R[\pi^{-1}]$ and $\F:=R/\pi R$.
Let $X_R\to \Spec R$ be a projective smooth scheme over $R$.
Put $X_K:=X_R\times_R K$ and $X_\F:=X_R\times_R \F$.
There is the {\it boundary map} 
\begin{equation}\label{boundary0}
\partial:\CH^2(X_K,1;\Z/n\Z)\lra \Pic(X_\F)/n,\quad n\geq0.
\end{equation}
Let us recall the definition.
We freely use the identifications \eqref{gersten3} and \eqref{lands}. 
Let $\sum (f,D)\in \ker (d_1\ot\Z/n)$ where $D$ is an irreducible divisor on
$X_K$ and $f$ is a rational function on $D$.
Let $D_R$ be the Zariski closure of $D$ in $X_R$.
Let $\widetilde{D}_R\to D_R$ be the normalization and
$j_D:\widetilde{D}_R\to D_R\hra X_R$ the composition.
The cycle $$Z=\sum j_{D*}\Div_{\widetilde{D}_R}(f)\in Z^2(X_R)/n=Z_1(X_R)/n$$ is supported
on $X_\F$ since $d_2(\sum(f,D))=0$ in $Z^2(X_K)/n$.
Thus one can consider it to be a divisor on $X_\F$ and hence it determines an element $[Z]$
of $\Pic(X_\F)$. We then define
\begin{equation}\label{boundary1}
\partial(\sum(f,D)):=[Z]\in \Pic(X_\F)/n .
\end{equation}
It is simple to show that \eqref{boundary1} is well-defined,
namely it annihilates the image of $d_2\ot\Z/n$.

\section{Surface containing a curve with nodes}\label{nodesect}
Let $K$ be a finite extension of $\Q_p$, $R$ the ring of integers and $\F$ the residue field.
For a scheme $V_R$ over $R$, we write $V_K:=V_R\times_R K$ and $V_\F:=V_R\times_R \F$.

\subsection{Conditions (A) and (B)}
Let $X_R\subset \P^3_R$ be a hypersurface which is smooth over $R$ and 
$C_R\subset X_R$ a hyperplane section which is flat over $R$.
Let $\pi_0:\tilde{C}_R\to C_R$ be the normalization
and $i:C'_R\to \tilde{C}_R$ a desingularization,
i.e. $C'_R$ is a regular arithmetic surface which is proper flat over $R$ and
$C'_K\os{\sim}{\to}\tilde{C}_K$.
Put $\pi:=\pi_0 i$.
$$
\xymatrix{
 C'_R \ar[r]^i \ar[rd]_\pi& \tilde{C}_R\ar[d]^{\pi_0}\\
 &C_R\ar[r]^\subset&X_R.
}
$$
Let
$$
C_{\F}=\sum_{j=1}^N D_j,\quad
\tilde{C}_{\F}=\sum_{j=1}^N \tilde{D}_j,\quad
C'_{\F}=\sum_{j=1}^N D'_j+\sum_l E_l
$$
be the irreducible decompositions such that $\pi_0(\tilde{D}_j)=D_j$,
$i(D'_j)=\tilde{D}_j$ and $E_l$ are exceptional curves (see Remark \ref{irred11}).
Note that $C_\F$ and $\tilde{C}_\F$ are reduced schemes (Remark \ref{reduced12}).

We consider the following two conditions on $(X_R,C_R)$.
\begin{description}
\item[(A)]
$X_K$ satisfies the following.

There is a nonsingular scheme $S$ over $\Q$ and a morphism $X_S\to S$ which has a Cartesian diagram
$$
\xymatrix{
&X_K\ar[d] \ar@{}[rd]|{\square}\ar[r]&X_S\ar[d]\\
&\Spec K\ar[r]&S.
}
$$
induced from an embedding $\Q(S)\hookrightarrow K$ such that
the complexes
\begin{equation}\label{sym1}
0\lra H^{i,2-i}\lra
H^{i-1,3-i}\ot \Omega^1_S
\quad (i=1,2)
\end{equation}
\begin{equation}\label{sym2}
H^{2,0}\lra
H^{1,1}\ot \Omega^1_S\lra
H^{0,2} \ot \Omega^2_S
\end{equation}
\begin{equation}\label{sym3}
0\lra
H^{2,0}\ot \Omega^1_S\lra
H^{1,1} \ot \Omega^2_S
\end{equation}
induced from the Gauss-Manin connection are exact at the middle terms.
Here we put $H^{i,j}
=(R^jf_*\Omega^i_{X_S/S})_{\mathrm{prim}}$
the Hodge $(i,j)$-component of the
primitive cohomology $H^2_\dR(X_S/S)_{\mathrm{prim}}:=H^2_\dR(X_S/S)/[H]$
with $H$ a hyperplane section.
\item[(B)] $C_R$ satisfies the following conditions {\bf (B-1)} and {\bf (B-2)}.
\item[(B-1)]
\begin{enumerate}
\item[(1)]$C_R$ is an irreducible (hence integral) scheme (cf. Remark \ref{reduced12}).
\item[(2)]$C_K$ has singular points $A_1,\cdots,A_m$
which are $K$-rational nodes.
\item[(3)]
$\pi^{-1}(A_i)$ consists of two $K$-rational points $P_i$ and $Q_i$.
\item[(4)]$C_{\F}$ is not irreducible
(hence so is neither $\tilde{C}_{\F}$ nor $C'_\F$, cf. Remark \ref{irred11}).
\end{enumerate}
\item[(B-2)]
There are $r_i\in \Z_p$ ($1\leq i \leq m$) which satisfy the following.
\begin{enumerate}
\item[(1)]
Let $J=J(C'_K)$ be the Jacobian variety of $C'_K$
and $\AJ:\CH_0(C'_K)_{\deg=0}\os{\sim}{\to} J(K)$
the Abel-Jacobi map. Then
\begin{equation}\label{4-2}
\sum_{i=1}^m r_i\AJ(P_i-Q_i)=0 \mbox{ in } \varprojlim_n J(K)/p^n.
\end{equation}
\item[(2)]
Let $P_{i,R}$ and $Q_{i,R}$ be the Zariski closure of $P_i$ and $Q_i$
in $C'_R$ respectively.
We denote by $(-\cdot -)_{C'_R}$ the 
intersection pairing on the arithmetic surface $C'_R$ (cf. \cite{lang} Chapter III).
Then there is no solution $x_l\in \Q_p$ which satisfy all of the following
equalities
\begin{equation}\label{4-1}
\begin{cases}
\sum_l x_l(E_l\cdot D'_k)_{C'_R}+\sum_{i=1}^m r_i(P_{i,R}-Q_{i,R}\cdot D'_k)_{C'_R}= 0& \forall k,\\
\sum_l x_l(E_l\cdot E_s)_{C'_R}+\sum_{i=1}^m r_i(P_{i,R}-Q_{i,R}\cdot E_s)_{C'_R}= 0& \forall s.
\end{cases}
\end{equation}
In other words, there are $q_k,q'_s\in \Z_p$ such that
\begin{equation}\label{4-1-1}
\sum_k q_k(E_l\cdot D'_k)_{C'_R}+\sum_s q'_s(E_l\cdot E_s)_{C'_R}=0,\quad \forall ~l
\end{equation}
\begin{equation}\label{4-1-2}
\sum_k q_k\left(\sum_{i=1}^m r_i(P_{i,R}-Q_{i,R}\cdot D'_k)_{C'_R}\right)
+\sum_{s} q'_s\left(\sum_{i=1}^m r_i(P_{i,R}-Q_{i,R}\cdot E_s)_{C'_R}\right)\ne 0.
\end{equation}
In case $C'_R=\tilde{C}_R$ (i.e. $\{E_l\}=\emptyset$), this simply means 
$$\sum_{i=1}^m r_i(P_{i,R}-Q_{i,R}\cdot D'_k)_{C'_R}\ne 0,\quad \exists ~k.$$
\end{enumerate}
\end{description}

\begin{rem}\label{intrem}
Since $P_{i,R}$ and $Q_{i,R}$ are $R$-sections,
they meet $C'_\F$ at nonsingular points transversally.
Thus we have
$$
(P_{i,R}\cdot D)_{C'_R}=
\begin{cases}
1& P_{i,R}\cap D\ne \emptyset\\
0& P_{i,R}\cap D= \emptyset\\
\end{cases},\quad
(Q_{i,R}\cdot D)_{C'_R}=
\begin{cases}
1& Q_{i,R}\cap D\ne \emptyset\\
0& Q_{i,R}\cap D= \emptyset\\
\end{cases}
$$
for any component $D\subset C'_\F$. 
\end{rem}

\begin{rem}\label{irred11}
There is 1-1 correspondence between irreducible components of $C_\F$ and those of $\tilde{C}_\F$
(hence $C_\F=\sum_j \pi_0(\tilde{D}_j)$ is the irreducible decomposition).
In fact, $\pi_0$ is a finite morphism (EGA IV 7.8) 
which is an isomorphism over the regular locus $C_R^{\mathrm{reg}}$.
Therefore $\pi^{-1}_0(C_\F^{\mathrm{reg}})\to C_\F^{\mathrm{reg}}$ is an isomorphism
and $\tilde{C}_\F-\pi_0^{-1}(C_\F^{\mathrm{reg}})$ is a finite set of closed points. 
This implies that there is a 1-1 correspondence between generic points of $C_\F$ and those of
$\tilde{C}_\F$.
\end{rem}
\begin{rem}\label{reduced12}
One can see that
$C_\F$ (and hence $\tilde{C}_\F$) are automatically reduced.
In fact let $(x,y,z,w)$ be homogeneous coordinates of $\P^3_\F$
such that $C_\F$ is defined by $w=0$. 
Write the defining equation of $X_\F$ by $F(x,y,z)+wG(x,y,z,w)$.
If $F$ has a decomposition $F=F_1^2F_2$, then the zero locus $\{F_1=G=w=0\}$
turns out to be a singular locus of $X_\F$. This contradicts with the assumption that
$X_\F$ is smooth.
\end{rem}
\begin{rem}\label{divisors}
The curves $\{D_j\}_j$ are linearly independent in $\Pic(X_{\ol{\F}})\ot\Q$.
To see it it is enough to show that the matrix
$M=((D_j\cdot D_k)_{X_\F})_{1\leq j,k\leq N}$ is nondegenerate
where $(-\cdot -)_{X_{\ol{\F}}}$ denotes the intersection pairing on $X_{\ol{\F}}$.
Let $e_j$ be the degree of $D_j$.
By definition,
$$%\begin{equation}\label{remark211}
e_j=(D_j\cdot C_\F)_{X_{\ol{\F}}}
=(D_j\cdot \sum_{k=1}^N D_k)_{X_{\ol{\F}}}
=D_j^2+\sum_{k\ne j}
(D_j\cdot D_k)_{X_{\ol{\F}}}.
$$%\end{equation}
We claim $(D_j\cdot D_k)_{X_{\ol{\F}}}=e_je_k$ for $j\ne k$.
In fact let $F(x,y,z)+wG(x,y,z,w)$ be the defining equation of $X_\F$ as before.
Let $F=\prod_{j=1}^N F_j$ be the irreducible decomposition.
Let $P\in D_j\cap D_k$. Since $X_\F$ is smooth, $G(P)\ne0$.
One has
$$
\O_{X_\F,P}/(F_j,F_k)\cong 
\O_{\P^3,P}/(F_j,F_k, F+wG)\cong 
\O_{\P^2,P}/(F_j,F_k). 
$$
Therefore $(D_j\cdot D_k)_{X_{\ol{\F}}}
=(D_j\cdot D_k)_{\P^2}=e_je_k$
by the theorem of B\'ezout.
Now we have 
$$
(D_j\cdot D_k)_{X_{\ol{\F}}}=\begin{cases}
e_j(1-\sum_{l\ne j}e_l)& j=k\\ 
e_je_k & j\ne k
\end{cases}
$$
and hence $$\det M=e_1\cdots e_N (1-\sum_j e_j)^{N-1}\ne0.$$
\end{rem}
\subsection{Infinite $p$-primary torsion in the Chow group of $0$-cycles}

\begin{thm}\label{main1}
Let $(X_R,C_R)$ satisfies {\bf (A)} and {\bf (B)}.
Then the corank of $\CH_0(X_K)[p^\infty]$ is nonzero (hence it is infinite).
\end{thm}
This follows from Bloch's exact sequence \eqref{bloch2} and 
the following Lemmas \ref{as1} and \ref{as2}.
\begin{lem}\label{as1}
The image of the regulator map $$
\CH^2(X_K,1)\ot\Q_p\lra H^1(G_K,H^2_\et(\ol{X},\Q_p(2)))
$$
coincides with the image of the map
$H^1(G_K,\Q_p(1))\to H^1(G_K,H^2_\et(\ol{X},\Q_p(2)))$ 
induced from the hyperplane class $\Q_p(1)\hra H^2_\et(\ol{X},\Q_p(2))$.
In other words the map
$$
\CH^2(X_K,1)\ot\Q_p\lra H^1(G_K,H^2_\et(\ol{X},\Q_p(2))_{\mathrm{prim}})
$$
is zero.
\end{lem}
\begin{pf}
The condition {\bf (A)} implies the above assertion (the proof is the same as in \cite{AS}).
\end{pf}
\begin{lem}\label{as2}
The corank of the image of the boundary map
$$
\partial:\CH^2(X_K,1;\Q_p/\Z_p)\to \Pic(X_\F)\ot\Q_p/\Z_p
$$
is greater than $1$.
In other words the image of the map
$$
\ol{\partial}: \CH^2(X_K,1;\Q_p/\Z_p)\lra (\Pic(X_\F)/[H])\ot\Q_p/\Z_p
$$
has nonzero corank where $H$ is a hyperplane section.
\end{lem}
\begin{pf}
The proof makes use of the condition {\bf (B)}.
We first construct an element $\xi_n\in \CH^2(X_K,1;\Z/p^n\Z)$.
By {\bf (B-2)} \eqref{4-2}, there is a rational function
$f_n$ on $C'_K$ such that
$$
\Div_{C'_K}(f_n)=\sum_{i=1}^mr_i (P_i-Q_i)\mod p^nZ_0(C'_K)
$$
where $Z_0(C'_K)$ denotes the free abelian group of closed points on
$C'_K$. Since
$$\pi_*\Div_{C'_K}(f_n)=\sum_{i=1}^m r_i(A_i-A_i)=0 \mod
p^n Z_0(X_K)$$ the pair $(f_n,C_K)$ determines an element $\xi_n\in
\CH^2(X_K,1;\Z/p^n\Z)$. By replacing $f_n$ with $cf_n$ for some constant $c\in K^*$, 
we may assume
that the support of $\Div_{C'_R}(f_n)$ does not containes 
the component $D'_1$ so that we have
$$
\Div_{C'_R}(f_n)\equiv \sum_{i=1}^m r_i(P_{i,R}-Q_{i,R})+\sum_{j>1} n_j D'_j+\sum_l m_l E_l 
\mod p^nZ_1(C'_R).
$$ 
By definition of the boundary map
$$
\partial(\xi_n)=\sum_{j>1}n_j[D_j]\quad \text{in }\Pic(X_\F)/p^n
$$
where $[D_j]$ denotes the cycle class of the divisor $D_j$ in $\Pic(X_\F)$
(cf. \S \ref{boundarysect}).
Note that $\{[D_j]\}_{j\geq1}$ are linearly independent in $\Pic(X_\F)\ot\Q$
(cf. Remark \ref{divisors}).
Therefore it is enough to show 
that $\min\{\ord_p(n_j)\}_j$ is bounded as $n\to+\infty$ (since 
$\ord_p(0):=+\infty$ by convension, it implies that some $n_j$ is nonzero).
Since $$(\Div_{C'_R}(f_n)\cdot D)_{C'_R}=0$$ for any components $D$ of $C'_\F$
and the intersection numbers on $C'_R$ are integers (\cite{lang} III \S 3) , one has
$$
\begin{cases}
\sum_j n_j(D'_j\cdot D'_k)_{C'_R}+\sum_l m_l(E_l\cdot D'_k)_{C'_R}+
\sum_{i=1}^m r_i(P_{i,R}-Q_{i,R}\cdot D'_k)_{C'_R}\equiv 0& \forall k,\\
\sum_j n_j(D'_j\cdot E_s)_{C'_R}+\sum_l m_l(E_l\cdot E_s)_{C'_R}+
\sum_{i=1}^m r_i(P_{i,R}-Q_{i,R}\cdot E_s)_{C'_R}\equiv 0& \forall s\\
\end{cases}
$$
modulo $p^n$.
It follows from {\bf (B-2)} \eqref{4-1-1} and \eqref{4-1-2} that we have
\begin{multline*}
-\sum_j n_j\left(\sum_{k} q_k(D'_j\cdot D'_k)_{C'_R}
+\sum_{s} q'_s(D'_j\cdot E_s)_{C'_R}\right)=\\
\sum_k q_k\left(\sum_{i=1}^m r_i(P_{i,R}-Q_{i,R}\cdot D'_k)_{C'_R}\right)
+\sum_{s} q'_s\left(\sum_{i=1}^m r_i(P_{i,R}-Q_{i,R}\cdot E_s)_{C'_R}\right)\ne 0
\end{multline*}
Here $n_j$ may depend on $n$. However so does none of $r_i$, $q_k$, $q'_s$ or
the above intersection numbers.
Therefore at least one $n_j$ is nonzero and its $p$-adic order is bounded
by that of the right hand side. This completes the proof.
\end{pf}

\begin{cor}\label{main1cor}
Let $(X_R,C_R)$ satisfies {\bf (A)} and {\bf (B)}.
Then the corank of $\CH_0(X_K\times_K L)[p^\infty]$ is nonzero for arbitrary finite extension $L/K$.
\end{cor}
\begin{pf}
It is enough to show that $X_L:=X_K\times_K L$ satisfies both of the assertions in Lemmas \ref{as1} and \ref{as2}.
Since {\bf (A)} is clearly satisfied for $X_L$, Lemma \ref{as1} holds.
The assertion in Lemma \ref{as2} remains true if we replace $X_K$ with $X_L$ because of the commutative
diagram
$$
\xymatrix{
\CH^2(X_K,1;\Q_p/\Z_p)\ar[r]^{\ol{\partial}}\ar[d]& (\Pic(X_\F)/[H])\ot\Q_p/\Z_p
\ar@{^{(}-_{>}}[d]\\
\CH^2(X_L,1;\Q_p/\Z_p)\ar[r]^{\ol{\partial}}& (\Pic(X_{\F_L})/[H])\ot\Q_p/\Z_p.
}
$$
Here the right vertical arrow is injective modulo finite group.
This completes the proof.
\end{pf}

\section{Construction of Quintic surface}\label{quinticsect}
In this section we construct $(X_{\Z_p},C_{\Z_p})$ which satisfies the conditions
{\bf (A)} and {\bf (B)} ($N=2$, $m=4$) in case $X_{\Z_p}$ is a quintic surface.

\subsection{Setting}\label{notationsect}
Let 
$$
t=(a_0,a_1,a_2,b_0,b_1,b_2,c_0,c_1,c_2, d_I)_{I=(i_0,i_1,i_2,i_3)}
$$
be the homogeneous coordinates of $\P_\Z^{43}$ where $I$ runs over the multi-indices
such that $i_k\geq0$ and $i_0+i_1+i_2+i_3=4$.
Put 
\begin{align*}
&H(x,y,z,w,t):=\sum_I d_I x^{i_0}y^{i_1}z^{i_2}w^{i_3},\\
&G(x,y,z,t):=x^2(x-z)^2L_1+y^2(y-z)^2L_2+
xy(x-z)(y-z)L_3
\end{align*}
where
$$
L_1=a_0x+a_1y+a_2z,\quad
L_2=b_0x+b_1y+b_2z,\quad
L_3=c_0x+c_1y+c_2z.
$$
We then consider a quintic homogeneous polynomial
\begin{equation}\label{noteq0}
F(x,y,z,w,t):=G(x,y,z,t)+w H(x,y,z,w,t)
\end{equation}
parametrized by $t$.
For an open set $S\subset \P^{43}_\Z$ we put
$$
X_S:=\{(x,y,z,w)\times t\in \P^3_\Z\times S~|~F(x,y,z,w,t)=0\}
$$
$$
C_S:=\{(x,y,z)\times t\in \P^2_\Z\times S~|~G(x,y,z,t)=0\}= X_S\cap \{w=0\}.
$$
We thus have a family of quintic surface containing a quintic curve which
has 4-nodes at $(x,y,z)=(0,0,1),(0,1,1),(1,0,1),(1,1,1)$:
\begin{equation}\label{family0}
\xymatrix{
 C_S \ar@{^{(}-_{>}}[rr] \ar[rd]&& X_S\ar[ld]\\
 &S.
}
\end{equation}
Hereafter we take $S$ to be an affine open set of $\P^{43}_\Z$
(which is of finite type over $\Z$) such that 
$X_S\to S$ is smooth.
Let 
$$
\begin{cases}
d_2(w)=-(a_0+a_2)w^2-(c_0+c_2)w-(b_0+b_2)\\
d_1(w)=a_0w^3-(a_1-c_0)w^2+(b_0-c_1)w-b_1\\
d_0(w)=(a_1+a_2)w^3+(c_1+c_2)w^2+(b_1+b_2)w
\end{cases}
$$
$$
\begin{cases}
e_2(u)=-(a_0+a_2)-(c_0+c_2)u-(b_0+b_2)u^2\\
e_1(u)=a_0-(a_1-c_0)u+(b_0-c_1)u^2-b_1u^3\\
e_0(u)=(a_1+a_2)u+(c_1+c_2)u^2+(b_1+b_2)u^3
\end{cases}
$$
and 
\begin{align*}
U_1&:=\{[s_0:s_1]\times w\in \P^1\times \bA^1~|~d_2(w)s_1^2+d_1(w)s_1s_0+d_0(w)s_0^2=0\},
\\
U_2&:=\{[t_0:t_1]\times u\in \P^1\times \bA^1~|~e_2(u)t_0^2+e_1(u)t_0t_1+e_0(u)t_1^2=0\}
\end{align*}
where $\P^1=\Proj\O(S)[x_0,x_1]$ and $\bA^1=\Spec\O(S)[z]$.
We glue $U_1$ and $U_2$ by identification
$$
\quad [t_0:t_1]\times u=[s_1:ws_0]\times w^{-1}
$$
and obtain a scheme $\tilde{C}_S$.
Put $s:=s_1/s_0$ and $t:=t_1/t_0$. Hereafter we simply denote the coordinates $[s_0:s_1]\times w$
and $[t_0:t_1]\times u$ by $(s,w)$ and $(t,w)$ respectively.
There is a finite morphism 
$\tilde{C}_S\to\P^1_S$ of degree 2 given by $(s,w)\mapsto w$.
The generic fiber of $\tilde{C}_S\to S$ is a 
nonsingular hyperelliptic curve of genus 2.
There is the normalization
$\pi_S:\tilde{C}_S\to C_S$ given by
$$
(s,w)\longmapsto (x,y,z)=
(s(w-s),w-s,w-s^2)
$$
$$
(t,u)\longmapsto (x,y,z)=
(t(1-t),u(1-t),u-t^2).
$$
Let $A_1=(0,0,1)$, $A_2=(0,1,1)$, $A_3=(1,0,1)$, $A_4=(1,1,1)$ be
the 4-nodes of $C_S$.
We put $\pi_S^{-1}(A_i)=\{P_i,Q_i\}$:
\begin{align}
&\begin{cases}
P_1(s,w)=(\alpha_1,\alpha_1)\\
Q_1(s,w)=(\alpha_2,\alpha_2)
\end{cases}
\quad
a_2\alpha_i^2+c_2\alpha_i+b_2=0,\label{PQ-1}\\
&\begin{cases}
P_2(s,w)=(0,\beta_1)\\
Q_2(s,w)=(0,\beta_2)
\end{cases}
\quad
(a_1+a_2)\beta_i^2+(c_1+c_2)\beta_i+b_1+b_2=0,\label{PQ-2}\\
&\begin{cases}
P_3(t,u)=(0,\gamma_1)\\
Q_3(t,u)=(0,\gamma_2)
\end{cases}
\quad
(b_0+b_2)\gamma_i^2+(c_0+c_2)\gamma_i+a_0+a_2=0,\label{PQ-3}\\
&\begin{cases}
P_4(s,w)=(1,\delta_1)\\
Q_4(s,w)=(1,\delta_2)
\end{cases}
\quad
(a_0+a_1+a_2)\delta_i^2+(c_0+c_1+c_2)\delta_i+b_0+b_1+b_2=0.\label{PQ-4}
\end{align}
We fix a regular affine scheme $T$ of finite type over $\Z$ and a generically finite morphism 
$T\to S$ such that
$P_i$ and $Q_i$ become $T$-valued points of $\tilde{C}_S$.
In other words, the function field $\Q(T)$ contains $\Q(S)$ and all of
$\alpha_i$,$\cdots$, $\delta_i$ in \eqref{PQ-1},$\cdots$,\eqref{PQ-4}.
Put $\tilde{C}_T:=\tilde{C}_S\times_S T$:
$$
\xymatrix{
&\tilde{C}_T\ar[d] \ar@{}[rd]|{\square}\ar[r]&\tilde{C}_S\ar[d]\\
&T\ar[r]\ar@/^/[u]^{P_i,Q_i}&S.
}
$$

\begin{thm}\label{main2}
There exists an embedding $\sigma:\O(T)\hra \Z_p$ such that 
the pair $(X_{\Z_p},C_{\Z_p})=(X_S\times_\sigma\Z_p,C_S\times_\sigma\Z_p)$
satisfies {\bf (A)} and {\bf (B)}.
\end{thm}
The rest of this section is devoted to prove Theorem \ref{main2}.
Combining the above with Corollary \ref{main1cor}, we have
\begin{thm}[Theorem \ref{main0}]\label{main00}
There exists a quintic surface $X_{\Q_p}\subset \P^3_{\Q_p}$ over $\Q_p$ such that
$\CH_0(X_{\Q_p}\times_{\Q_p} K)[p^\infty]$ is infinite for arbitrary finite extension $K/\Q_p$.
\end{thm}
\subsection{Condition (A)}\label{Asect}

\begin{prop}\label{symmetrizer}
Let $X_S/S$ be as above.
Then by shrinking $S_\Q:=S\times_\Z\Q$ to a small open set if necessary
the sequences \eqref{sym1}, \eqref{sym2} and \eqref{sym3} are exact.
Thus for any embedding $\Q(S)\hra K$, $X_K:=X_S\times_S K$ satisfies the condition {\bf (A)}.
\end{prop}
\begin{pf}
Let 
$$
R_S=\O(S_\Q)[x,y,z,w]/(\frac{\partial F}{\partial x},\frac{\partial F}{\partial y},
\frac{\partial F}{\partial z},\frac{\partial F}{\partial w})
$$
be the Jacobian ring of the quintic polynomial $F$ \eqref{noteq0}.
It follows from the theory of Jacobian rings that one has
$$
R_S^{11-5i}\cong H^{i,2-i}:=F^{i}/F^{i+1}H_\dR(X_S/S)_{\mathrm{prim}}\ot\Q
$$
for $0\leq i \leq 2$. Moreover the tangent space of $S$ is canonically isomorphic to
the homogeneous part $I^5$ of degree 5 of ideal 
$$
I=\langle x^2(x-z)^2,y^2(y-z)^2,xy(x-z)(y-z),w\rangle
$$
of $R_S$
and the Gauss-Manin connection can be identified with the dual of the ring product
$R_S^\bullet\ot I^5\to R_S^{\bullet+5}$.
Thus to show the exactness of \eqref{sym1}, \eqref{sym2} and \eqref{sym3} 
it is enough to show the following 4-sequences are exact:
\begin{equation}\label{jac1}
R_S^{1+5i}\ot I^5\lra R_S^{6+5i}\lra 0\quad (i=0,1)
\end{equation}
\begin{equation}\label{jac3}
R^1_S\ot \bigwedge^2 I^5\lra R^6_S\ot I^5 \lra R^{11}_S
\end{equation}
\begin{equation}\label{jac4}
R^6_S\ot \bigwedge^2 I^5\lra R^{11}_S\ot I^5 \lra 0
\end{equation}
All of them can be checked by direct calculation (with the aid of computer).
\end{pf}

\subsection{Condition (B)}\label{Bsect}
Let $\sigma:\O(T)\to \Q_p$ be a ring homomorphism.
Set $X^\sigma_{\Z_p}:=X_S\times_\sigma \Z_p$,
$X^\sigma_{\Q_p}:=X_S\times_\sigma \Q_p$,
$X^\sigma_{\F_p}:=X_S\times_\sigma \F_p$ and similarly for $C_S$ and $\tilde{C}_S$.
Put $P_{i,\Z_p}^\sigma:=P_i\times_\sigma \Z_p$, 
$Q_{i,\Z_p}^\sigma:=Q_i\times_\sigma \Z_p$ and
$P_{i,\Q_p}^\sigma:=P_i\times_\sigma \Q_p$, 
$Q_{i,\Q_p}^\sigma:=Q_i\times_\sigma \Q_p$.

\bigskip

We consider the following conditions on $\sigma$.
Write $a^\sigma=\sigma(a)$.
\begin{description}
\item[(i)]
The homomorphism $\sigma$ is injective. In other words,
$(a_1/a_0)^\sigma,\cdots,(d_I/a_0)^\sigma$ are algebraically independent over $\Q$.
\item[(ii)]
$X^\sigma_{\F_p}$ is smooth over $\F_p$ (hence $X^\sigma_{\Z_p}$ is smooth over $\Z_p$).
\item[(iii)]
$\tilde{C}_{\F_p}^\sigma$ has two irreducible components $D_1$ and $D_2$.
\item[(iv)]
$\tilde{C}_{\Z_p}^\sigma$ is an integral regular scheme (hence $\tilde{C}_{\Q_p}$ is a nonsingular hyperelliptic curve
of genus 2).
\item[(v)]
There is a subset $\{i_1,i_2,i_3\}\subset\{1,2,3,4\}$ such that
\begin{equation}\label{iv-1}
\begin{cases}
P_{i,\Z_p}^\sigma\cap D_1\ne \emptyset\\
P_{i,\Z_p}^\sigma\cap D_2= \emptyset\\
\end{cases}
\quad\mbox{and}\quad
\begin{cases}
Q_{i,\Z_p}^\sigma\cap D_1= \emptyset\\
Q_{i,\Z_p}^\sigma\cap D_2\ne \emptyset\\
\end{cases}
\end{equation}
or
\begin{equation}\label{iv-2}
\begin{cases}
P_{i,\Z_p}^\sigma\cap D_1= \emptyset\\
P_{i,\Z_p}^\sigma\cap D_2\ne \emptyset\\
\end{cases}
\quad\mbox{and}\quad
\begin{cases}
Q_{i,\Z_p}^\sigma\cap D_1\ne \emptyset\\
Q_{i,\Z_p}^\sigma\cap D_2= \emptyset\\
\end{cases}
\end{equation}
for any $i\in \{i_1,i_2,i_3\}$.
\item[(vi)]
Let $J^\sigma=J(\tilde{C}_{\Q_p}^\sigma)$ be the Jacobian variety and
$\AJ:\CH_0(\tilde{C}^\sigma_{\Q_p})_{\deg=0}\os{\sim}{\to} J^\sigma(\Q_p)$ the Abel-Jacobi map.
Let $\{i_1, i_2, i_3\}$ be as in {\bf (vi)}.
We may assume that \eqref{iv-1} holds by exchanging $P^\sigma_{i,\Z_p}$ and $Q^\sigma_{i,\Z_p}$ if necessary.
Then $\AJ(P^\sigma_{i_1,\Q_p}-Q^\sigma_{i_1,\Q_p})-\AJ(P^\sigma_{i_3,\Q_p}-P^\sigma_{i_3,\Q_p})$ and 
$\AJ(P^\sigma_{i_2,\Q_p}-Q^\sigma_{i_2,\Q_p})-\AJ(P^\sigma_{i_3,\Q_p}-P^\sigma_{i_3,\Q_p})$ are linearly independent
over $\Q_p$ in $J^\sigma(\Q_p)\ot\Q\cong\Q_p^2$.
\end{description}

\begin{prop}\label{exists}
There exists $\sigma$ which satisfies all of the conditions {\bf (i)},$\cdots$,{\bf (vi)}.
\end{prop}
We shall give a proof of Proposition \ref{exists} in \S \ref{existssect}.

\begin{prop}\label{exists1}
Suppose that $\sigma$ satisfies all of the conditions {\bf (i)},$\cdots$,{\bf (vi)}.
Then $(X_{\Z_p},C_{\Z_p})=(X_{\Z_p}^\sigma,C_{\Z_p}^\sigma)$ satisfies the condition {\bf (B)}.
\end{prop}
\begin{pf}
{\bf (B-1)} is straightforward.
We see {\bf (B-2)}.
Since $J$ is a 2-dimensional abelian variety over $\Q_p$, one has
$$
J^\sigma(\Q_p)\cong \Z_p^2+\mbox{(finite group)}
$$
by the theorem of Mattuck \cite{mattuck} (see also \cite{serre} Part II Ch.V \S 7, Corollary 4).
Therefore there is a nontrivial relation
\begin{equation}\label{mattuck}
\sum_{k=1}^3 r_k\AJ(P^\sigma_{i_k,\Q_p}-Q^\sigma_{i_k,\Q_p})=0\quad \mbox{in } \plim n J^\sigma(\Q_p)/p^n
\quad (r_i\in \Z_p)
\end{equation}
which gives the condition {\bf (B-2)} (1).
We show that $r_1$, $r_2$ and $r_3$ satisfy {\bf (B-2)} (2).
Since $\tilde{C}^\sigma_{\Z_p}$ is regular, what we want to show is that
\begin{equation}
\sum_{j=1}^3 r_i(P^\sigma_{i_k,\Z_p}-Q^\sigma_{i_k,\Z_p}\cdot D_l)_{\tilde{C}^\sigma_{\Z_p}}
=\pm\sum_{j=1}^3r_i
\end{equation} 
is nonzero
(see Remark \ref{intrem} for the above equality).
Suppose $r_1+r_2+r_3=0$.
It follows from \eqref{mattuck}
that $\AJ(P^\sigma_{i_1,\Q_p}-Q^\sigma_{i_1,\Q_p})-\AJ(P^\sigma_{i_3,\Q_p}-Q^\sigma_{i_3,\Q_p})$
and $\AJ(P^\sigma_{i_2,\Q_p}-Q^\sigma_{i_2,\Q_p})-\AJ(P^\sigma_{i_3,\Q_p}-Q^\sigma_{i_3,\Q_p})$ are 
not linearly independent in $J^\sigma(\Q_p)\ot\Q\cong \Q_p^2$.
This contradicts with {\bf (vi)}.
\end{pf}
Theorem \ref{main2} follows from Propositions \ref{symmetrizer}, \ref{exists}
and \ref{exists1}.

\subsection{Proof of Proposition \ref{exists}}\label{existssect}
For a smooth scheme $V$ over $\Q_p$, we denote by $V^\an$ a topological space $V(\Q_p)$
endowed with the $p$-adic manifold structure
(cf. \cite{serre} Part II Chapter III).
For a smooth scheme $V$ over $k\subset \Q_p$, we simply write $V^\an=(V\times_k\Q_p)^\an$.

We put 
\begin{align*}
U_2&=\{\sigma:\O(T)\to \Z_p~|~\sigma\mbox{ satisfies {\bf (ii)}}\}\subset T^\an \\
&\vdots\\
U_5&=\{\sigma:\O(T)\to \Z_p~|~\sigma\mbox{ satisfies {\bf (v)}}\}\subset T^\an,
\end{align*}
and
\begin{align*}
F_1&=\{\sigma:\O(T)\to \Q_p~|~\sigma\mbox{ satisfies {\bf (i)}}\}\subset T^\an \\
U_6&=\{\sigma:\O(T)\to \Q_p~|~\sigma\mbox{ satisfies {\bf (vi)}}\}\subset T^\an .
\end{align*}
Our goal is to show $F_1\cap U_2\cap \cdots \cap U_6\ne\emptyset.$
\begin{lem}\label{eslem1}
$F_1$ is a dense subset of $T^\an$.
Namely for any open ball ${\mathbb B}$ in $T^\an$, one has ${\mathbb B}\cap F_1\ne \emptyset$.
\end{lem}
\begin{pf}
Easy.
\end{pf}
For ring homomorphisms $\sigma,\tau:\O(T)\to\Z_p$,
we say $\sigma\equiv \tau$ mod $p^m$ if 
$x^\sigma\equiv x^{\tau}$ mod $p^n$
for all $x\in \O(T)$.
\begin{lem}\label{exisistslem0}
\begin{enumerate}
\item[(1)] Suppose $\sigma\equiv \tau$ mod $p$.
Then $\sigma$ satisfies {\bf (ii)} (resp. {\bf (iii)}, {\bf (v)}) if and only if so does $\tau$.
\item[(2)] Suppose $\sigma\equiv \tau$ mod $p^2$.
Then $\sigma$ satisfies {\bf (iv)} if and only if so does $\tau$.
\end{enumerate}
Thus $U_2$, $U_3$, $U_4$ and $U_5$ are open sets of $T^\an$.
\end{lem}
\begin{pf}
(1) is clear.
We see (2). It is easy to see that $\tilde{C}^\sigma_{\Z_p}$ is connected.
It is enough to see whether $\tilde{C}^\sigma_{\Z_p}$ is regular or not around singular points of 
$\tilde{C}^\sigma_{\F_p}$.
Let $x$ be a singular point of $\tilde{C}^\sigma_{\F_p}$.
Let $f_i$ be the local equations of the irreducible component $D_i$ of $\tilde{C}^\sigma_{\F_p}$ around $x$.
Let $F_i$ be a lifting of $f_i$ to characteristic zero.
Then one can write the local equation of $\tilde{C}^\sigma_{\Z_p}$ as
$
F_1F_2+pG
$.
$\tilde{C}_{\Z_p}^\sigma$ is regular around $x$ if and only if $G(x)\not\equiv 0$ mod $p$.
Thus the assertion follows.
\end{pf}
It is not difficult to construct $\sigma:\O(S)\to \Z_p$
which satisifies {\bf (ii)}, $\cdots$, {\bf (v)}.
Hence
\begin{lem}\label{eslem2}
$U_2\cap U_3\cap U_4 \cap U_5\ne \emptyset$.
\end{lem}
\begin{lem}\label{eslem3}
$U_6$ is a dense open set of $T^\an$.
\end{lem}

To prove this, we prepare some notations.
Let $T_\Q:=T\times_\Z\Q$ and $\tilde{C}_{S_\Q}:=\tilde{C}_S\times_\Z\Q$.
Let $\M_2$ be the moduli scheme of curves of genus 2 over $\Q$, and
$\uC\to \M_2$ the universal curve.
Then there is a dominant morphism $S_\Q\to \M_2$ such that
$\tilde{C}_{S_\Q}\cong S_\Q\times_{\M_2}\uC$. 
Note that $T_\Q$ is a 8-dimensional variety and $\M_2$ is a 3-dimensional variety.
Put $\tilde{C}_{T_\Q}:=\tilde{C}_{S}\times_{S} T_\Q$:
$$
\xymatrix{
\tilde{C}_{T_\Q}\ar[d]\ar[r]\ar@{}[rd]|{\square}&\tilde{C}_{S_\Q}
\ar[r]\ar[d]\ar@{}[rd]|{\square}& \uC\ar[d]\\
T_\Q\ar[r]&S_\Q\ar[r]& \M_2.
}
$$
Let $\uJ\to \M_2$ be the Jacobian of $\uC$.
The divisor $P_i-Q_i$ induces the morphism
 \begin{equation}\label{fi-1}
f_i:T_\Q\lra \uJ.
\end{equation}

\medskip

We first prove that $U_6$ is an open set.
Recall the theorem of Mattuck.
Let $\uG$ be the Lie algebra bundle of $\uJ$ over $\M_2$.
We endow $\uG$ with the $p$-adic topology, and denote it by $\uG^\an$. 
Then there is a subbundle $\Lambda\subset\uG^\an$
whose fiber is isomorphic to
$\Z_p^2$ and a subgroup bundle $G\subset \uJ^\an$ of finite index such that
 \begin{equation}\label{gi-5}
\xymatrix{
\Lambda \ar[r]^{\exp}_{\sim} \ar[d]_\cap& G\ar[d]^{\cap}\\
\uG^\an&\uJ^\an
}
\end{equation}
where ``$\exp$" is the exponential map (\cite{serre} Part II Ch.V \S 7).
In particular there is a $p$-adically continuous homomorphism
 \begin{equation}\label{fi-2}
\varepsilon:\uJ^\an\lra \Lambda
\end{equation}
whose kernel and cokernel are finite. Clearly it is an open map on $p$-adic manifolds.
Put $g_1:=f_{i_1}-f_{i_3}$, $g_2:=f_{i_2}-f_{i_3}$ and
$g:=g_1\times g_2:T_\Q\to \uJ\times \uJ
$.
We denote by $g^\an_i$ etc. the associated $p$-adic analytic map:
\begin{equation}\label{gi-3}
g_i^\an:T^\an\lra \uJ^\an,\quad
g^\an:T^\an\lra \uJ^\an\times \uJ^\an.
\end{equation}
Letting \begin{equation}\label{gi-V}
V=\{(v_1,v_2)\in \Lambda\times \Lambda~|~v_1\wedge v_2\ne0\}
\end{equation}
be a dense open set,
we have
$$
U_6=(g^\an)^{-1}(\varepsilon \times \varepsilon )^{-1}(V).
$$
This shows that $U_6$ is an open set.

\medskip

Next we show that $U_6$ is a dense subset.
The map $g$ gives rise to the map
\begin{equation}\label{gi-2}
(g_*)_x:\mathrm{tan}(T_\Q)_x\lra \bigoplus_{i=1}^2 \mathrm{tan}(\uJ)_{g_i(x)}
\end{equation}
of the Zariski tangent space at a point $x\in T_\Q$.

\begin{claim}\label{surjectivity}
There is a closed subset $Z\subsetneq T_\Q$ such that
\eqref{gi-2} is surjective for all $x\not\in Z$.
\end{claim}
\begin{pf}
Changing the variable $s$ with $Y$ by $Y=2d_2(w)s+d_1(w)$
one has the Weierstrass form of the hyperelliptic curve $C_{T_\Q}$:
$$
Y^2=f(w)=d_1(w)^2-4d_0(w)d_2(w)=a_0^2w^6+\cdots+b_1^2.
$$
We want to show that the pull-back
\begin{equation}\label{gi-1-p}
\bigoplus_{i=1}^2 g_i^*\Omega^1_{\uJ/\M_2} \lra \Omega^1_{T_\Q/\M_2},\quad
(\omega_1,\omega_2)\longmapsto g_1^*\omega_1+g_2^*\omega_2
\end{equation}
of K\"ahler differentials is injective at the generic point of $T_\Q$.
Evaluating $a_0=0$, one has a closed subscheme $T_0\hra T_\Q$ and
\begin{equation}\label{igusa}
C_{T_\Q}\times_{T_\Q} T_0: Y^2=f_0(w)=v_0w^5-v_1w^4+\cdots -v_5,
\end{equation}
$$
\begin{cases}
v_0=4a_2(a_1+a_2)\\
v_1=-((a_1+c_0)^2+4a_2(c_0+c_1+c_2)+4c_2(a_1+a_2))\\
\qquad\vdots
\end{cases}
$$
Still $T_0\to \M_2$ is dominant.
We show that the composition map
\begin{equation}\label{gi-1}
\bigoplus_{i=1}^2 g_i^*\Omega^1_{\uJ/\M_2}|_{T_0} \lra \Omega^1_{T_\Q/\M_2}|_{T_0}\lra
 \Omega^1_{T_0/\M_2}
\end{equation}
is bijective at the generic point of $T_0$, which implies the injectivity of \eqref{gi-1-p} and 
hence the desired assertion.
Note that $\Omega^1_{\uJ/\M_2}$ is a locally free sheaf of rank 2 generated by
invariant 1-forms $\frac{dw}{Y}$ and $w\frac{dw}{Y}$. One has
\begin{align*}
f_1^*\frac{dw}{Y}&=
\frac{d\alpha_1}{2\alpha_1d_2(\alpha_1)+d_1(\alpha_1)}
-\frac{d\alpha_2}{2\alpha_2d_2(\alpha_2)+d_1(\alpha_2)}\\
f_1^*w\frac{dw}{Y}&=
\frac{\alpha_1d\alpha_1}{2\alpha_1d_2(\alpha_1)+d_1(\alpha_1)}
-\frac{\alpha_2d\alpha_2}{2\alpha_2d_2(\alpha_2)+d_1(\alpha_2)}\\
f_2^*\frac{dw}{Y}&=
\frac{d\beta_1}{d_1(\beta_1)}
-\frac{d\beta_2}{d_1(\beta_2)}\\
f_2^*w\frac{dw}{Y}&=
\frac{\beta_1d\beta_1}{d_1(\beta_1)}
-\frac{\beta_2d\beta_2}{d_1(\beta_2)}\\
f_3^*\frac{dw}{Y}&=
\frac{\gamma_1d\gamma_1}{e_1(\gamma_1)}
-\frac{\gamma_2d\gamma_2}{e_1(\gamma_2)}\\
f_3^*w\frac{dw}{Y}&=
\frac{d\gamma_1}{e_1(\gamma_1)}
-\frac{d\gamma_2}{e_1(\gamma_2)}\\
f_4^*\frac{dw}{Y}&=
\frac{d\delta_1}{2d_2(\delta_1)+d_1(\delta_1)}
-\frac{d\delta_2}{2d_2(\delta_2)+d_1(\delta_2)}\\
f_4^*w\frac{dw}{Y}&=
\frac{\delta_1d\delta_1}{2d_2(\delta_1)+d_1(\delta_1)}
-\frac{\delta_2d\delta_2}{2d_2(\delta_2)+d_1(\delta_2)}.
\end{align*}
We want to show that
\begin{equation}\label{igusa-1}
f_{i_1}^*\frac{dw}{Y}-f_{i_3}^*\frac{dw}{Y},\quad
f_{i_1}^*w\frac{dw}{Y}-f_{i_3}^*w\frac{dw}{Y},\quad
f_{i_2}^*\frac{dw}{Y}-f_{i_3}^*\frac{dw}{Y},\quad
f_{i_2}^*w\frac{dw}{Y}-f_{i_3}^*w\frac{dw}{Y}
\end{equation}
generate $\Omega^1_{T_0/\M_2}$ as $\O_{T_0}$-module at the generic point.

The affine coordinate ring of $\M_2$ is described by
Igusa's $j$-invariants $J_2$, $J_4$, $J_6$, $J_8$ and $J_{10}$ (\cite{igusa}).
In particular, $\Omega^1_{\M_2/\Q}$ is generated by
\begin{equation}\label{igusa-2}
d(J_4/J_2^2),\quad d(J_6/J_2^3),\quad d(J_{10}/J_2^5)
\end{equation}
generically.
Therefore it is enough to show that \eqref{igusa-1} and \eqref{igusa-2}
generate $\Omega^1_{T_0/\Q}$ as $\O_{T_0}$-module at the generic point.
We know the explicit forms of $J_2$, $\cdots$, $J_{10}$ (see \S \ref{app} Appendix).
Therefore one can check it by direct calculations
(the details are left to the reader since they are long and tedious).
\end{pf}
We prove that $U_6$ is dense in $T^\an$.
Let $\sigma\in T^\an$ be an arbitrary point.
For any open ball $\mathbb B$ about $\sigma$, we want to show ${\mathbb B}\cap U_6\ne \emptyset$.
There is a point $\sigma_0\in {\mathbb B}-Z(\Q_p)$.
So it is enough to show ${\mathbb B}_0\cap U_6\ne\emptyset$ for any open ball ${\mathbb B}_0$ about $\sigma_0$. 
It follows from Claim \ref{surjectivity}
that the map $g^\an$ \eqref{gi-3} is an open map on (sufficiently small) ${\mathbb B}_0$. In particular 
$g^\an({\mathbb B}_0)$ is an open subset of $\uJ^\an\times \uJ^\an$ and hence 
$(\varepsilon\times\varepsilon)g^\an({\mathbb B}_0)$ is also open in $\Lambda\times \Lambda$.
Since $V$ \eqref{gi-V} is a dense open set, we have $V\cap (\varepsilon\times\varepsilon)g^\an({\mathbb B}_0)\ne\emptyset$. 
This means
${\mathbb B}_0\cap U_6={\mathbb B}_0\cap (g^\an)^{-1}(\varepsilon\times\varepsilon)^{-1}(V)\ne\emptyset$,
which is the desired assertion.
This completes the proof of Lemma \ref{eslem2}.

\bigskip

We finish the proof of Proposition \ref{exists}.
It follows from Lemmas \ref{eslem2} and \ref{eslem3} that
$U_2\cap\cdots\cap U_6$ is a nonempty open set.
By Lemma \ref{eslem1} one has $F_1\cap U_2\cap\cdots\cap U_6 \ne \emptyset$.
This is the desired assertion.

\section{Appendix: Igusa's $j$-invariants}\label{app}

In \cite{igusa}, Igusa gave the arithmetic invariants of hyperelliptic curves of genus 2.
For a hyperelliptic curve which has an affine equation
$$
y^2=f(x)=v_0x^5-v_1x^4+v_2x^3-v_3x^2+v_4x-v_5
$$
they are given as follows.
\begin{align*}
J_2=&5v_0v_4-2v_1v_3+4^{-1}3v_2^2\\
J_4=&-8^{-1}[25v_0^2v_3v_5-15v_0^2v_4^2-15v_0v_1v_2v_5+7v_0v_1v_3v_4
+2^{-1}v_0v_2^2v_4\\&-v_0v_2v_3^2+4v_1^3v_5-v_1^2v_2v_4-v_1^2v_3^2
+v_1v_2^2v_3-3\cdot 2^{-4}v_2^4]\\
J_6=&-16^{-1}[2^{-1}5^3v_0^3v_2v_5^2-25v_0^3v_3v_4v_5+5v_0^3v_4^3-25v_0^2v_1^2v_5^2\\
&-10v_0^2v_1v_2v_4v_5+10v_0^2v_1v_3^2v_5-v_0^2v_1v_3v_4^2
-4^{-1}5v_0^2v_2^2v_3v_5\\&-4^{-1}11v_0^2v_2^2v_4^2+2^{-1}7v_0^2v_2v_3^2v_4
-v_0^2v_3^4+6v_0v_1^3v_4v_5\\&-3v_0v_1^2v_2v_3v_5+2^{-1}7v_0v_1^2v_2v_4^2
-2v_0v_1^2v_3^2v_4+3\cdot 4^{-1}v_0v_1v_2^3v_5\\&-4^{-1}7v_0v_1v_2^2v_3v_4
+v_0v_1v_2v_3^3+7\cdot 16^{-1}v_0v_2^4v_4-4^{-1}v_0v_2^3v_3^2-v_1^4v_4^2\\
&+v_1^3v_2v_3v_4-4^{-1}v_1^2v_2^3v_4
-4^{-1}v_1^2v_2^2v_3^2+8^{-1}v_1v_2^4v_3-2^{-6}v_2^6]\\
J_8=&4^{-1}[J_2J_6-J_4^2]\\
J_{10}=&v_0D
\end{align*}
where
$$
D=\begin{vmatrix}
v_0&-v_1&v_2&-v_3&v_4&-v_5&0&0&0\\
0&v_0&-v_1&v_2&-v_3&v_4&-v_5&0&0\\
0&0&v_0&-v_1&v_2&-v_3&v_4&-v_5&0\\
0&0&0&v_0&-v_1&v_2&-v_3&v_4&-v_5\\
5v_0&-4v_1&3v_2&-2v_3&v_4&0&0&0&0\\
0&5v_0&-4v_1&3v_2&-2v_3&v_4&0&0&0\\
0&0&5v_0&-4v_1&3v_2&-2v_3&v_4&0&0\\
0&0&0&5v_0&-4v_1&3v_2&-2v_3&v_4&0\\
0&0&0&0&5v_0&-4v_1&3v_2&-2v_3&v_4\\
\end{vmatrix}
$$
is Sylvester's resultant.

Putting the degree of $J_{2i}$ to be $2i$,
the affine coordinate ring of $\M_2$ is given by the homogeneous part of degree $0$
in the graded ring $\Q[J_2,\cdots,J_{10},J_{10}^{-1}]$ (\cite{igusa} Theorem 2).

\noindent
Department of Mathematics, Hokkaido University,
Sapporo 060-0810,
JAPAN

\smallskip

\noindent
{\it E-mail} : \textbf{asakura@math.sci.hokudai.ac.jp}

\end{document}